# Boundary Behavior and Confinement of Screw Dislocations

Marco Morandotti[1]
[1]Fakultät für Mathematik, Technische Universität München, Boltzmannstrasse, 3,
85748 Garching bei München, Germany.

**ABSTRACT**

In this note we discuss two aspects of screw dislocations dynamics: their behavior near the boundary and a way to confine them inside the material. In the former case, we obtain analytical results on the estimates of collision times (one dislocation with the boundary and two dislocations with opposite Burgers vectors with each other); numerical evidence is also provided. In the latter, we obtain analytical results stating that, under imposing a certain type of boundary conditions, it is energetically favorable for dislocations to remain confined inside the domain.

**INTRODUCTION**

Introduced in a seminal paper by Volterra in 1907 [12], dislocations have been proposed as the mechanism that is ultimately responsible for plasticity in metals: Orowan [9], Polanyi [10], and Taylor [11] reached this conclusion independently in 1934. Yet, it was not until 1956 that dislocations were observed with the help of an electron microscope [5]. Nowadays, much effort is put into studying their behavior, and especially their dynamics.

We build on a model proposed in 1999 by Cermelli & Gurtin [3] for screw dislocations undergoing antiplane shear to study two different situations: the behavior of dislocations near the domain boundary and a situation in which dislocations can be confined inside the domain by imposing a suitable boundary condition. The theoretical framework adopted is that of linearized elasticity in the plane.

In both cases, the assumptions on the deformation allow to reduce the problem from a fully 3D one in the cylinder $\Omega \times \mathbb{R}$ to a 2D one in the cross-section $\Omega \subseteq \mathbb{R}^2$, which we will assume to be an open connected set with $C^2$ boundary that does not touch itself.

To study the behavior near the boundary, we provide accurate estimates on the Green's function for the laplacian in two dimensions, which will translate into estimates on the direction of the Peach-Koehler force acting on a dislocation near the boundary. The first result that we obtain (see Theorem 1) is that, if one dislocation is sufficiently close to the boundary and sufficiently far away from the others in the system, then its Peach-Koehler force is directed along the outward unit normal to the boundary at the closest point.

Next, we turn to estimates on the collision times: we prove that in the situation just described the dislocation which is closest to the boundary will collide with it in finite time. We can also provide estimates on the collision time of two dislocations with opposite Burgers vectors, provided that they are sufficiently close to each other and that the remaining dislocations are sufficiently far away. Furthermore, in both cases, if the other dislocations are sufficiently diluted, no other collision events will happen (see Theorems 2 and 3).

To study the confinement, we resort to variational methods and we prove that, under suitably chosen boundary condition, an energy functional that depends on the position of the dislocations is minimized when the dislocations are in the interior of the domain. Together with the *core radius approach* [1,3], the main tool used to prove the results is Γ-convergence (the results are in fact formulated in terms of continuous convergence, which is a stronger notion) [4]. We present here the statements of the results in the case of a single dislocation in the domain (see Theorem 4 and Corollary 5) and refer the reader to [8] for the case of many dislocations.

The results presented here have been obtained in collaboration with T. Hudson [7] and I. Lucardesi, R. Scala, and D. Zucco [8].

**BOUNDARY BEHAVIOR**

The *renormalized energy* (see, *e.g.*, [1,3]) for a system of $n$ dislocations $z_1,..., z_n \in \Omega$ subject to *natural boundary conditions* may be expressed as

$$\mathcal{E}_n(z_1,..., z_n) := \sum_{i<j} b_i b_j \left( k_\Omega(z_i, z_j) - \frac{1}{2\pi} \log|z_i - z_j| \right) + \frac{1}{2} \sum_{i=1}^n h_\Omega(z_i),$$

where $k_\Omega(x, y)$ is related to the Green's function of the laplacian in the plane via $G_\Omega(x, y) = -\frac{1}{2\pi} \log|x - y| + k_\Omega(x, y)$ and solves

$$\begin{cases} -\Delta_x k_\Omega(x, y) = 0, & \text{in } \Omega, \\ k_\Omega(x, y) = \frac{1}{2\pi} \log |x - y|, & \text{on } \partial\Omega, \end{cases}$$

and $h_\Omega(x) := k_\Omega(x, x)$. The function $h_\Omega$ solves $-\Delta_x h_\Omega(x) = \frac{2}{\pi} e^{-4\pi h_\Omega(x)}$ and it has been studied in [2]. The Peach-Koehler force acting on a dislocation [6] is obtained by taking the negative of the gradient of the renormalized energy with respect to the dislocation position: $f_i(z_1,..., z_n) = -\nabla_{z_i} \mathcal{E}_n(z_1,..., z_n)$, for $i = 1,..., n$, and it is responsible for the motion of dislocations via the law

$$\dot{z}_i(t) = f_i(z_1,..., z_n), \text{ for } i = 1,..., n,$$

complemented with an initial condition at time $t = 0$.

The assumptions on the regularity of $\Omega$ imply that it satisfies a *uniform interior and exterior disk condition* with radius $\bar{\rho}$. Define the function $d_n : \Omega^n \to [0, +\infty)$ as

$$d_n(x_1,..., x_n) := \begin{cases} \text{dist}(x_1, \partial\Omega), & n = 1, \\ \min_i \text{dist}(x_i, \partial\Omega) \wedge \min_{i \neq j} |x_i - x_j|, & \text{otherwise}, \end{cases}$$

and, given $0 < \delta < \gamma < \text{diam}\Omega/2$, let

$$\mathcal{D}_{n,\delta,\gamma} := \{(z_1, z') \in \Omega^n | d_1(z_1) < \delta, d_{n-1}(z') > \gamma\}.$$

Estimates on $\nabla h_\Omega$, derived from estimates on $\nabla G_\Omega$, yield the following results.
**Theorem 1** (free boundaries attract dislocations [7, Theorem 3.1]). *Let $n \in \mathbb{N}$, let $\sigma \in (0,1)$, and*

let $\delta \in (0, \sigma\bar{\rho})$ and $\gamma \in (\max\{2\delta, \bar{\rho}\}, \text{diam}\Omega/2)$. Let $z = (z_1, z') \in \mathcal{D}_{n,\delta,\gamma}$. Then, if $s \in \partial\Omega$ is the boundary point closest to $z_1$, the Peach-Koehler force $f_1(z)$ on the dislocation $z_1$ satisfies

$$f_1(z) = \frac{\nu(s)}{4\pi d_1(z_1)} + O(1),$$

where $\nu(s)$ is the outward unit normal to $\partial\Omega$ at $s$ and $O(1)$ denotes a term which is uniformly bounded for all $z \in \mathcal{D}_{n,\delta,\gamma}$.

Theorem 1 is proved by obtaining an estimate of the type

$$\left| f_1(z) - \frac{\nu(s)}{4\pi d_1(z_1)} \right| \leq \frac{C_{n,\sigma}(\gamma)}{2\pi\bar{\rho}},$$

where $C_{n,\sigma}(\gamma)$ is a constant depending on the geometry of the domain $\Omega$ via $\bar{\rho}$ and on how far all the other dislocations are from $z_1$ and from $\partial\Omega$.

**Theorem 2** (collision with the boundary [7, Theorem 3.2]). *Let $n \in \mathbb{N}$, let $\sigma \in (0,1)$, and let $\gamma_0 > 0$. There exists $\delta_0 > 0$ such that, if $z(0) \in \mathcal{D}_{n,\delta_0,\gamma_0}$, then there exists $T_{coll} > 0$ such that the evolution $z(t)$ is defined for $t \in [0, T_{coll}]$, $z(t) \in \Omega^n$ for $t \in [0, T_{coll})$, and $z_1(T_{coll}) \in \partial\Omega$ and $z'(T_{coll}) \in \Omega^{n-1}$.*
*Moreover, the following estimate on the collision time holds: $T_{\text{coll}} \leq 2\pi\delta_0^2 + O(\delta_0^3)$.*

Let $\zeta, \eta > 0$ with $\zeta < \eta$ and define

$$\mathcal{C}_{n,\zeta,\eta} := \{(z_1, z_2, z'') \in \Omega^n \mid |z_1 - z_2| < \zeta, d_{n-2}(x'') > \eta, \text{dist}(\{z_1, z_2\}, \{z_3, \ldots, z_n\} \cup \partial\Omega) > \eta\}.$$

**Theorem 3** (collision between dislocations [7, Theorem 3.4]). *Let $n \in \mathbb{N}$ $(n \geq 2)$ and let $\eta_0 \in (0, \text{diam}\Omega/2)$. There exists $\zeta_0 > 0$ such that, if $z(0) \in \mathcal{C}_{n,\zeta_0,\eta_0}$, then there exists $T'_{coll} > 0$ such that the evolution $z(t)$ is defined for $t \in [0, T'_{coll}]$, $z(t) \in \Omega^n$ for $t \in [0, T'_{coll})$, and $z_1(T'_{coll}) = z_2(T'_{coll}) \in \Omega$ and $z''(T'_{coll}) \in \Omega^{n-2}$.*
*Moreover, the following estimate on the collision time holds: $T'_{\text{coll}} \leq \frac{\pi\zeta_0^2\eta_0^2}{2(\eta_0^2 - \zeta_0^2 - 2(n-2)\zeta_0\delta_0)}$.*

**Analytically solvable cases and numerical results**

We consider different domains where computations can be carried out analytically, namely the unit disk, the half plane, and the plane. The last two are not admissible for our theorems, but a close inspection of the renormalized energy shows that the evolution is well defined also in the case of unbounded domains. Explicit computations can be done for: one dislocation $z$ in (i) the half plane and (ii) the unit disk ($\bar{\rho} = 1$), (iii) two dislocations $z_1, z_2$ of opposite Burgers modulus $b_1 = +1 = -b_2$ in the unit disk ($\bar{\rho} = 1$), and (iv) two dislocations in the plane. In cases (i) and (ii) the dislocation hits the boundary in finite time; case (iii) shows a variety of scenarios, i.e., collisions with the boundary, collision between dislocations, and unstable equilibria; in case (iv) the dislocations will collide, or the evolution exists for all time. We refer the reader to [7, Section 4] for details.

We include plots from numerical simulations of the dynamics in different scenarios. Figure 1(a) shows the superposition of 5000 runs of the scenario of case (iii), where initial conditions $(z_1(0), z_2(0))$ have been randomly generated in $\mathcal{D}_{2,0.2,0.5}$. Figure 1(b) shows an histogram of hitting times, which agree, at leading order, with the bound provided by Theorem 2 $T_{\text{coll}} \leq 2\pi\delta_0^2 \approx 0.2513$. Figure 1(c) shows plots of 80 trajectories of one dislocation evolving in the cardioid, which has an unstable equilibrium point at its center: the initial conditions are chosen on a circle of radius 0.1 centered at the equilibrium point. Due to the interaction with the boundary, the dislocation starts following a curved line and then hits the boundary perpendicularly (up to numerical artifacts), as indicated in Theorem 1.

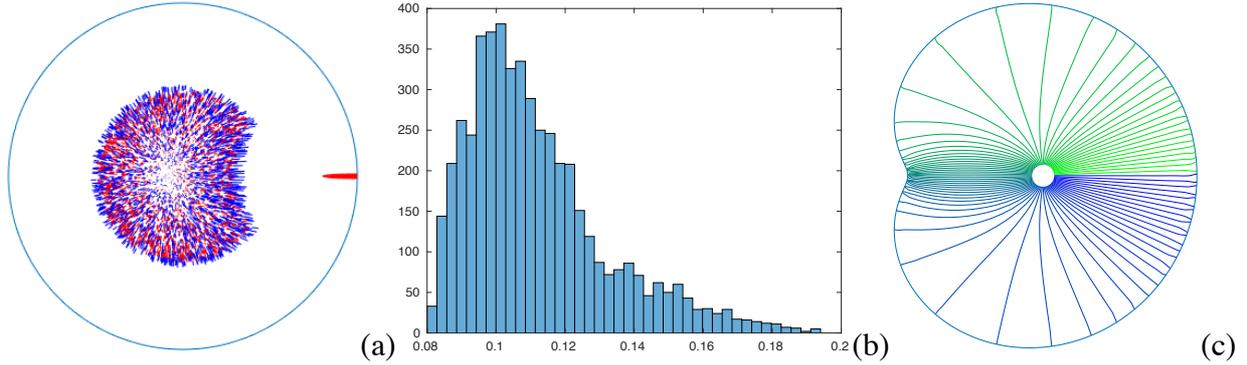

**Figure 1**. (a) superposition of 5000 runs for two dislocations in the unit disk: red corresponds to Burgers modulus $b_i = +1$, blue to $b_i = -1$; (b) histogram of hitting times; (c) cardioid.

**CONFINEMENT**

In order to state the results about the confinement of screw dislocations inside a crystal, we have to take the point of view of the energy (as a function of the position of the dislocations) and show that it is minimized if the dislocations are inside the domain, provided that a suitable boundary condition is imposed on $\partial\Omega$. To this aim, we prescribe an external tangential strain on $\partial\Omega$, by means of a function $f \in L^1(\partial\Omega)$ such that $\int_{\partial\Omega} f(x) \, d\mathcal{H}^1(x) = 2\pi$. This choice will allow for at most one dislocation to form in the crystal; we denote by $a \in \Omega$ the position of the dislocation, so that the strain of the deformed crystal is represented by a field $F_a \in L^1(\Omega; \mathbb{R}^2) \cap L^2_{\text{loc}}(\Omega \setminus \{a\}; \mathbb{R}^2)$, solution to

$$\begin{cases} \text{div} F_a = 0, & \text{in } \Omega, \\ \text{curl} F_a = 2\pi\delta_a, & \text{in } \Omega, \\ F_a \cdot \tau = f, & \text{on } \partial\Omega, \end{cases}$$

where $\delta_a$ is the Dirac delta centered at $a$ and $\tau$ is the tangent unit vector to $\partial\Omega$. The energy associated with such a system is the linearized elastic one, which, in terms of $F_a$, is written $\frac{1}{2}\int_\Omega |F_a|^2 \, dx$. It is easy to show that this energy diverges to $+\infty$, but by means of the core radius approach, the problem can be regularized. To this end, consider $\epsilon > 0$ and denote by $B_\epsilon(a)$ the disk of radius $\epsilon$ centered at $a$, by $\bar{B}_\epsilon(a)$ its closure, and by $\Omega_\epsilon(a) := \Omega \setminus \bar{B}_\epsilon(a)$. Now, $F_a$ can be approximated in $L^2_{\text{loc}}(\Omega \setminus \{a\}; \mathbb{R}^2)$ by a sequence of fields $F_a^\epsilon \in L^2(\Omega_\epsilon(a); \mathbb{R}^2)$ which solve

$$\begin{cases} \text{div} F_a^\epsilon = 0, & \text{in } \Omega_\epsilon(a), \\ \text{curl} F_a^\epsilon = 0, & \text{in } \Omega_\epsilon(a), \\ F_a^\epsilon \cdot \tau = f, & \text{on } \partial\Omega \setminus \bar{B}_\epsilon(a), \\ F_a^\epsilon \cdot \nu = 0, & \text{on } \partial B_\epsilon(a) \cap \Omega, \end{cases}$$

where $\nu$ is the outer unit normal to $\partial\Omega$. The system above characterizes the minimizers of the energy functional

$$\mathcal{E}_\epsilon(a) := \min\left\{\frac{1}{2}\int_{\Omega_\epsilon(a)} |F|^2 \, dx : F \in L^2(\Omega_\epsilon(a); \mathbb{R}^2), \text{curl} F = 0, F \cdot \tau = f \text{ on } \partial\Omega \setminus \bar{B}_\epsilon(a)\right\}.$$

The functionals of which we study the $\Gamma$-convergence are $\mathcal{F}_\epsilon(a) := \mathcal{E}_\epsilon(a) - \pi|\log \epsilon|$.

**Theorem 4** ([8, Theorem 1.1]). *The functionals $\mathcal{F}_\epsilon$ continuously converge in $\bar{\Omega}$ to the functional $\mathcal{F}: \bar{\Omega} \to (-\infty, +\infty]$ defined, for $a \in \Omega$, as*

$$\mathcal{F}(a) := \pi \log d_1(a) + \frac{1}{2}\int_{\Omega_{d_1(a)}(a)} |K_a + \nabla v_a|^2 dx + \frac{1}{2}\int_{B_{d_1(a)}(a)} |\nabla v_a|^2 dx,$$

*and $\mathcal{F}(a) = +\infty$ otherwise. Here $K_a(x) := \rho_a^{-1}(x)\hat{\theta}_a(x)$ (with $(\rho_a, \theta_a)$ polar coordinates centered at $a$) and $v_a$ is the solution to*

$$\begin{cases} \Delta v_a = 0, & \text{in } \Omega, \\ v_a = g - \theta_a, & \text{on } \partial\Omega, \end{cases}$$

*for $g$ a primitive of the boundary datum $f$. In particular, $\mathcal{F}$ is continuous over $\bar{\Omega}$ and diverges to $+\infty$ as the dislocation approaches $\partial\Omega$, that is $\mathcal{F}(a) \to +\infty$ as $d(a) \to 0$. As a consequence, $\mathcal{F}$ attains its minimum in the interior of $\Omega$.*

**Corollary 5** ([8, Corollary 1.2]). *There exists $\epsilon_1 > 0$ such that, for every $\epsilon \in (0, \epsilon_1)$, the problem*

$$\inf\{\mathcal{E}_\epsilon(a) : a \in \bar{\Omega}\}$$

*admits a minimizer only in the interior of $\Omega$. Moreover, if $a^\epsilon \in \Omega$ is a minimizer, then (up to subsequences) $a^\epsilon \to a$ and $\mathcal{F}_\epsilon(a^\epsilon) \to \mathcal{F}(a)$ as $\epsilon \to 0$, where $a$ is a minimizer of the $\Gamma$-limit $\mathcal{F}$ defined in Theorem 4. In particular, for $\epsilon$ small enough, all the minimizers stay uniformly (with respect to $\epsilon$) away from the boundary.*

Theorem 4 and Corollary 5 express that confinement of dislocations can be obtained by imposing a traction condition at the boundary. The $\Gamma$-limit $\mathcal{F}$ of the regularized functionals $\mathcal{E}_\epsilon$ is a functional that attains its minimum in the interior of the domain $\Omega$. Moreover, in Corollary 5 it is stated that the minimizers of the energy stay well separated from $\partial\Omega$, guaranteeing that the dislocation will not collide with the boundary.

Theorem 4 and Corollary 5 have their counterpart for $n$ dislocation $a_1, \ldots, a_n \in \Omega$. The statements are more technical and are far beyond the scope of this note. We refer the reader to [8, Theorem 1.3 and Corollary 1.4] for a precise statement and to [8, Section 5.4.3] for the proofs.

**DISCUSSION**

All the results presented here are valid under geometric assumptions on the domain that assure the boundedness of the curvature and the uniqueness of the boundary point of minimal distance for an interior point sufficiently close to the boundary. The regularity of $\Omega$ required in the introduction can be weakened, especially for the results concerning the confinement; the details can be found in the papers [7,8].

The results regarding the confinement rely on imposing the boundary condition $F_a \cdot \tau = f$. The mechanical meaning of this condition is that of imposing a sort of infinitesimal rotation at the boundary. The author and his collaborators are investigating if this is suitable for designing a real experimental setup.

**ACKNOWLEDGMENTS**

The author is grateful to SISSA, the Mathematics Institute of the University of Warwick, the Department of Mathematics of Politecnico di Torino, and the Faculty of Mathematics of Technische Universität München, where this research was conducted. The author's research was partially funded by the ERC Advanced grant *Quasistatic and Dynamic Evolution Problems in Plasticity and Fracture* (grant no.: 290888), by the ERC Starting grant *High-Dimensional Sparse Optimal Control* (grant no.: 306274), and by the 2015 GNAMPA project *Fenomeni critici nella meccanica dei materiali: un approccio variazionale*. The author is a member of the Gruppo Nazionale per l'Analisi Matematica, la Probabilità e le loro Applicazioni (GNAMPA) of the Istituto Nazionale di Alta Matematica (INdAM).